\pdfoutput=1
\RequirePackage{fix-cm}
\documentclass[12pt]{article}   

\usepackage{a4wide}
\usepackage{graphicx}

\usepackage[utf8]{inputenc}
\usepackage[english]{babel}
\usepackage{listings}   
\usepackage{algorithm2e}
\SetAlFnt{\footnotesize}
\SetAlCapFnt{\small}
\SetAlCapNameFnt{\small}

\usepackage{color}
\usepackage {amsmath,amsfonts,amssymb}
\usepackage{subfig}
\usepackage{texnames}
\usepackage{url}
\definecolor{hotpink}{rgb}{0.9,0,0.5}
\usepackage[colorlinks,urlcolor=blue,citecolor=hotpink,linkcolor=blue]{hyperref}

\newcommand{\F}{\mathbb{F}}
\newcommand{\fl}{\mathrm{fl}}
\newcommand{\R}{\mathbb{R}}

\renewcommand{\leq}{\leqslant}
\renewcommand{\geq}{\geqslant}
\renewcommand{\dim}{d}

\newcommand{\half}{\textstyle \frac12}

\begin{document}

\title{Efficient implementation of symplectic implicit Runge-Kutta schemes with simplified Newton iterations
}

\author{Mikel Anto\~nana, Joseba Makazaga,  Ander Murua \\
KZAA saila, Informatika Fakultatea, UPV/EHU\\ Donostia / San Sebasti\'an}

\date{  }

\maketitle

\begin{abstract}
We are concerned with the efficient implementation of symplectic  implicit Runge-Kutta (IRK) methods applied to systems of (non-necessarily Hamiltonian) ordinary differential equations by means of Newton-like iterations. We pay particular attention to symmetric symplectic IRK schemes (such as collocation methods with Gaussian nodes). 
For a $s$-stage IRK scheme used to integrate a $\dim$-dimensional system of ordinary differential equations, the application of simplified versions of Newton iterations requires solving at each step several linear systems (one per iteration) with the same  $s\dim \times s\dim$ real coefficient matrix.  We propose rewriting such $s\dim$-dimensional linear systems as an equivalent $(s+1)\dim$-dimensional systems that can be solved by performing  the LU decompositions of $[s/2] +1$ real matrices of size $\dim \times \dim$. 
We present a C implementation (based on Newton-like iterations) of Runge-Kutta collocation methods with Gaussian nodes  that make use of such a rewriting of the linear system and that takes special care in reducing the effect of round-off errors. We report some numerical experiments that demonstrate the reduced round-off error propagation of our implementation.
\end{abstract}

\section{Introduction}
\label{sec:intro}

The main goal of the present work is the efficient implementation of symplectic implicit Runge-Kutta (IRK) schemes for stiff ODE (non-necessarily Hamiltonian) problems. Our primary interest is on geometric numerical integration, which motivates us to solve  the implicit equations determining each step to full machine precision. The stiff character of the target problems leads us to solving the implicit equations by some modified version of Newton method. This typically requires repeatedly solving linear systems of equations with coefficient matrix of the form
\begin{equation}
\label{eq:Matrix}
  \left( I_s \otimes I_\dim  - h \, A \otimes J \right) \in \mathbb{R}^{s\dim \times s\dim}
\end{equation}
where $A\in \R^{s\times s}$ is the coefficient matrix of the RK scheme, and $J$ is some common approximation of the Jacobian matrices evaluated at the stage values. 

 An standard approach, independently introduced in~\cite{Liniger1970}, \cite{Butcher1976}, and~\cite{Bickart1977}, to efficiently solve such linear systems takes advantage of the special structure of the matrix (\ref{eq:Matrix}). More specifically, (\ref{eq:Matrix}) is similar to a block-diagonal matrix with $s$ blocks  $I_\dim - h\, \lambda_j J \in \R^{\dim \times \dim}$ ($j=1,\ldots,s$), one per eigenvalue $\lambda_j$ of $A$. Typically, the coefficient matrix $A$ of standard high order implicit RK schemes have $[s/2]$ complex conjugate pairs of eigenvalues (plus a real one for odd $s$).

The main contribution of the present paper is a technique for transforming $s\,\dim$-dimensional systems with coefficient matrix (\ref{eq:Matrix}) into an equivalent $(s+1)\dim$-dimensional systems that can be solved by performing the LU decomposition of $[s/2]+1$ real matrices of size $\dim \times \dim$ (and some additional multiplication of matrices of the same size). We pay particular attention to implicit Runge-Kutta schemes that are both symmetric and symplectic. However, our technique is also applicable for non-symmetric symplectic IRK schemes, and also for some symmetric non-symplectic IRK schemes (see last paragraph in Subsection~3.3 for more details).

In addition, we present an algorithm that implements symmetric sympletic IRK schemes  (such as RK collocation methods with Gaussian nodes)  by making use of the above technique. Special care is taken to try to reduce the effect of round-off errors by adapting some techniques used (for the implementation of sympletic IRK schemes with fixed-point iterations) in \cite{Antonana2017}. 

The plan of the paper is as follows: Section~2 summarizes some standard material about implicit Runge-Kutta methods and Newton-like iterations and fix some notation. Section~3 presents  our new technique to solve the simplified linear system of Newton iterations for symplectic IRK schemes. Section~4 is devoted to describe our implementation of symplectic IRK methods with Newton-like iterations. Some numerical results are reported in Section~5. A few concluding remarks can be found in Section~6.

\section{Implementation of Implicit Runge-Kutta schemes with Newton-like iterations}

\subsection{Implicit Runge-Kutta schemes}

We consider initial value problems of systems of ODEs of the form
\begin{equation}
\label{eq:ivp}
\frac{d}{dt}y=f(t,y),\quad  y(t_0)=y_0,
\end{equation}
where $f: \R^{\dim+1}\to \R^\dim$ is a sufficiently smooth map and $y_0 \in \R^\dim$.

Given a time discretization $t_0 < t_1 < t_2 < \cdots$, the numerical approximations  $y_{n} \approx y(t_{n})$, ($n=1,2,\ldots$) to the solution $y(t)$ of the initial value problem (\ref{eq:ivp}) is obtained by means of a one-step integrator as
\begin{equation}
\label{eq:one_step}
y_{n+1} = \Phi(y_{n},t_{n},t_{n+1}-t_{n}), \quad n=0,1,2,\ldots,
\end{equation}
for a map $\Phi: \R^{\dim+2} \to \R^\dim$ determined in some way from $f: \R^{\dim+1} \to \R^\dim$. 

In the case of a s-stage implicit Runge-Kutta method, the map $\Phi$ is determined in terms 
of the real coefficients $a_{ij}$ ($1 \leq i, j \leq s$) and $b_{i}$, $c_{i}$ ($1\leq i \leq s$)
as
\begin{equation}
\label{eq:PhiIRK}
\Phi(y,t,h) := y+ h\sum^s_{i=1} b_i \, f(t + c_i h, Y_{i}),
\end{equation}   
where the {\em stage vectors} $Y_{i}$ are implicitly defined as functions of $(y,t,h) \in \R^{\dim+2}$ by
\begin{equation}
\label{eq:Y}
Y_{i} =y+ h \sum^s_{j=1}{a_{ij}\,f(t + c_j h,Y_{j})}, \quad  i=1 ,\ldots, s.
\end{equation}
Typically, 
\begin{equation*}
c_i = \sum_{j=1}^{s} a_{i j}, \quad i=1,\ldots,s.
\end{equation*}

The equations (\ref{eq:Y}) can be solved for the stage vectors $Y_i$ by means of some iterative procedure, starting, for instance, from  $Y_{i}^{[0]}=y$, $i=1,2,\ldots,s$. ( In the non-stiff case, it is usually more efficient initializing the stage vectors with some other procedure that uses the stage values of the previous steps~\cite{Hairer2006}).

A very simple iterative procedure is fixed point iteration. For stiff problems, fixed-point iteration is not appropriate, and Newton iteration may be used to compute the stage vectors $Y_i$ from (\ref{eq:Y}).  For non-stiff problems, Newton iteration may still be an attractive option in some cases, in particular for very high precision computations (for quadruple precision or in arbitrary precision arithmetic calculations) if implemented with mixed-precision strategies~\cite{Baboulin2009} (which reduce the cost of the linear algebra and the evaluation of the Jacobians, performed in lower precision arithmetic than the evaluations of the right-hand side of the system of ODEs).

In any case, since at each Newton iteration
$s$ evaluations of the Jacobian matrix $\frac{\partial f}{\partial y}$ and a LU decomposition of a $s\dim \times s\dim$ matrix are required,  some computationally cheaper variants are often
 used instead.
 
\subsection{Newton-like iterations}

Recall that a Newton iteration can be used to compute for $k=1,2,\ldots$  the approximations $Y_i^{[k]}$ of $Y_i$ ($i=1,\ldots,s$) in  (\ref{eq:Y}) as follows: 
\begin{align}
\label{eq:(1)Newton_iteration}
1) & \quad r_i^{[k]} := -Y_{i}^{[k-1]} + y + h \sum_{j=1}^{s}\, a_{ij}\, f(t + c_j h,Y_{j}^{[k-1]}), \quad  i=1 ,\ldots, s, \\
\label{eq:(2)Newton_iteration}
\begin{split}
2) & \quad \mathrm{Solve \ } \Delta Y_{i}^{[k]} \mbox{\  from \ } \\
& \quad \Delta Y_{i}^{[k]}  - h \sum_{j=1}^{s}\, a_{ij}\, J_j^{[k]} \Delta Y_{j}^{[k]} = r_i^{[k]} \quad  i=1 ,\ldots, s, \\
& \mbox{where} \quad  J_i^{[k]}=\frac{\partial f}{\partial y}(t + c_i h,Y_{i}^{[k]}) \quad	\mbox{for} \quad  i=1,\ldots, s, 
\end{split} \\
\label{eq:(3)Newton_iteration}
3)& \quad Y_i^{[k]} := Y_i^{[k-1]} + \Delta Y_i^{[k]}, \quad  i=1 ,\ldots, s,
\end{align}
Observe that,  $s$ evaluations of the Jacobian matrix $\frac{\partial f}{\partial y}$ and a LU decomposition of a $s\dim \times s\dim$ matrix are required (in addition to $s$ evaluations of $f$) at each iteration. This is typically computationally too expensive, and some variants of the 
 full Newton algorithm are implemented instead. Among others, the following alternatives are possible:
 \begin{itemize}
\item Application of simplified Newton iterations. This consists on replacing the Jacobian matrices $J_i^{[k]}$ in (\ref{eq:(2)Newton_iteration}) by 
$J_i^{[0]}=\frac{\partial f}{\partial y}(t + c_i h,Y_{i}^{[0]})$. In that case, LU decomposition is done only once and the linear system 
\begin{equation}
\label{eq:lssn}
\Delta Y_{i}^{[k]}  - h \sum_{j=1}^{s}\, a_{ij}\, J_j^{[0]} \Delta Y_{j}^{[k]} = r_i^{[k]} \quad  i=1 ,\ldots, s,
\end{equation}
has to be solved at each of the simplified Newton iterations.  If the simple initialization $Y_{i}^{[0]}=y$ ($i=1,\ldots,s$) is considered (this is typically the case when solving stiff systems) and $f$ does not depend on $t$, then $J_i^{[0]}=J:= \frac{\partial f}{\partial y}(y)$ for each $i=1,\ldots,s$, and the linear system (\ref{eq:lssn})
reduces to 
\begin{equation}
\label{eq:lsssn}
 \left( I_s \otimes I_\dim  - h \, A \otimes J \right) \Delta Y^{[k]} = r^{[k]},
\end{equation}
where
\begin{equation*}
Y^{[k]} = \left(
\begin{matrix}
Y_1^{[k]}\\
\vdots\\
Y_s^{[k]}
\end{matrix}
\right) \in \R^{s\dim}, \quad 
r^{[k]} = \left(
\begin{matrix}
r_1^{[k]}\\
\vdots\\
r_s^{[k]}
\end{matrix}
\right)\in \R^{s\dim}.
\end{equation*}
Even in the case where some  initialization procedure other than $Y_{i}^{[0]}=y$ is used, in practice the linear system (\ref{eq:lssn}) is often replaced by (\ref{eq:lsssn}), where $J$ is some common approximation of $\frac{\partial f}{\partial y}(t + c_i h,Y_{i}^{[0]})$, $i=1,\ldots,s$. An appropriate choice~\cite{Xie2009} is
 $J:=  \frac{\partial f}{\partial y}(t+\bar c \, h,\bar y)$, where $\bar c = \frac{1}{s} \sum_{i=1}^{s}c_i$ (which for symmetric methods gives $\bar c = \frac12$) and $\bar y =  \frac{1}{s} \sum_{i=1}^{s}Y_i^{[0]}$.  
Often, it will be sufficient to evaluate instead of $\frac{\partial f}{\partial y}$ a computationally cheaper approximation of it. 

\item Applying the original Newton iteration by solving the linear systems (\ref{eq:(2)Newton_iteration}) with some iterative method~\cite{Saad2003} preconditioned by the inverse of the matrix 
\begin{equation}
\label{eq:matrix}
I_s \otimes I_\dim  - h \, A \otimes J.
\end{equation}
In practice,  the linear systems (\ref{eq:(2)Newton_iteration}) are only approximately solved with the iterative method. In such case, the resulting scheme is sometimes referred to as inexact Newton iteration~\cite{Saad2003}. Further variants of Newton-like iterations will be obtained if the Jacobian matrices are not updated at each iteration.
\end{itemize}
In any of the two alternatives above, one needs to repeatedly solve 
linear systems of the form
\begin{equation}
\label{eq:linsys}
 \left( I_s \otimes I_\dim  - h \, A \otimes J \right) \Delta Y = r,
\end{equation}
for given  $r \in \R^{s\dim}$. From now on, we will refer to (\ref{eq:linsys}) as simplified linear system (of Newton-like iterations).

Of course, (\ref{eq:linsys}) may be solved by previously computing the LU decomposition of the full $sd \times sd$ matrix (\ref{eq:matrix}), but this may be done more efficiently. 

An standard approach~\cite{Liniger1970}\cite{Butcher1976}\cite{Bickart1977} consists on diagonalizing the matrix $A$ as $\Lambda = S^{-1} A S=\mathrm{diag}(\lambda_1,\ldots,\lambda_s)$,  and computing the LU decomposition of the matrix
\begin{equation*}
I_s \otimes I_\dim  - h \, \Lambda \otimes J = (S^{-1} \otimes I_\dim) \left( I_s \otimes I_\dim  - h \, A \otimes J\right) (S \otimes I_\dim).
\end{equation*}
In that case, one needs to compute the LU decomposition of a real (resp. complex) $\dim \times \dim$ matrix for each distinct real eigenvalue (for each distinct pair of complex eigenvalues) of $A$.

 Alternatively, some authors~\cite{Brugnano2014}\cite{Jay2009} propose solving (\ref{eq:linsys}) by an iterative procedure preconditioned by the inverse of
\begin{equation}
\label{eq:matrix2}
I_s \otimes I_\dim  - h \, \tilde A \otimes J, 
\end{equation}
where  $\tilde A \in \R^{s \times s}$ is a matrix chosen so that the  LU decomposition of (\ref{eq:matrix2}) can be more efficiently computed than that of (\ref{eq:matrix}). 

In next section, we propose a new technique to efficiently solve simplified linear systems (\ref{eq:linsys}) of Newton iterations, provided  that the IRK scheme is symplectic.

\section{Efficient solution of simplified linear systems for symplectic IRK schemes}

\subsection{Symplectic IRK schemes}

In what follows, we consider symplectic IRK schemes,  that is~\cite{JMSanz-Serna1994}, IRK schemes whose coefficients satisfy
\begin{equation} \label{eq:sympl_cond_1}
b_{i}a_{ij}+b_{j}a_{ji}-b_{i}b_{j}=0, \ \ 1 \leq i,j \leq s.
\end{equation}

Condition (\ref{eq:sympl_cond_1}) guarantees that the discrete flow resulting from the application of the IRK scheme to an autonomous Hamiltonian system is symplectic, with important favorable consequences in the long-term behavior of the numerical solution~\cite{Hairer2006}. Condition (\ref{eq:sympl_cond_1}) also implies that, when applied to an ODE system with a quadratic invariant, then it is also a conserved quantity for the numerical solution provided by the IRK scheme.

Nevertheless, our interest in condition (\ref{eq:sympl_cond_1}) is of a completely different nature: We will see that such a condition allows to solve efficiently linear systems of the form (\ref{eq:linsys})
for a given $\dim \times \dim$ real matrix $J$ and a given  $r \in \R^{s\dim}$.
We will pay particular attention to symplectic IRK schemes that additionally satisfy (possibly after some reordering of the stage values $Y_i$) the symmetry condition~~\cite{Hairer2006}
\begin{equation}
\label{eq:symmetry_cond}
\begin{split}
b_{s+1-i} &= b_{i}, \quad   c_{s+1-i} = 1 - c_{i}, \quad 1 \leq i \leq s,\\
b_{j} &= a_{s+1-i,s+1-j} + a_{i, j}, \quad 1 \leq i,j \leq s,
\end{split}
\end{equation}
 In particular, the IRK schemes of collocation type with Gaussian nodes are both symplectic and symmetric.

\subsection{Alternative symplecticity and symmetry characterizations}

The map $\Phi$ determining the steps (\ref{eq:one_step}) of a IRK scheme can be alternatively written as
\begin{equation*}
\Phi(y,t,h) := y+ z,
\end{equation*}   
where the stage vectors $Y_{i} \in \R^\dim$ and the increment $z\in \R^\dim$ are implicitly defined as functions of $(y,t,h) \in \R^{\dim+2}$ by
\begin{align}
\label{eq:Y2}
Y_{i} &=y+ \frac{z}{2} + h \sum^s_{j=1}{\bar a_{ij}\,f(t + c_j h,Y_{j})}, \quad  i=1 ,\ldots, s, \\
\label{eq:z}
z &=  h\sum^s_{i=1} b_i \, f(t + c_i h, Y_{i}),
\end{align}
where
 \begin{equation}
\bar a_{i j} = a_{i j} - \frac{b_j}{2}, \quad 1 \leq i, j \leq s.
\end{equation}

Condition (\ref{eq:sympl_cond_1}) may be equivalently characterized 
in terms of the matrix $\bar A = (\bar a_{i j})_{i,j=1}^{s}$ and the diagonal matrix $B$ with diagonal entries $b_1,\ldots,b_s$.
Indeed, (\ref{eq:sympl_cond_1})  is equivalent to the requirement that the real $s\times s$ matrix $(B \bar A)$ be antisymmetric.

As for the symmetry condition (\ref{eq:symmetry_cond}), it reads
 \begin{equation}
\label{eq:symmetry_cond2}
\begin{split}
b_{s+1-i} &= b_{i}, \quad \bar  c_{s+1-i} = - \bar c_{i}, \quad 1 \leq i \leq s,\\
\bar a_{s+1-i,s+1-j} &= -\bar a_{i, j}, \quad 1 \leq i,j \leq s,
\end{split}
\end{equation}
where $\bar c_i = c_i -\frac12$ for $i=1,\ldots,s$.

\subsection{Efficient solution of the linear systems of the form (\ref{eq:linsys})}
\label{ss:sIRK}

From now on, we will only consider, without loss of generality\footnote{Any symplectic IRK method with $b_i=0$ for some $i$ is equivalent to a symplectic IRK scheme with fewer stages and $b_i\neq 0$ for all $i$~\cite{Hairer1994}}, symplectic IRK schemes with invertible $B$. Since $B \bar A$ is antisymmetric,  so is $B^{\frac12} \bar A B^{-\frac12}$, which implies that $\bar A$ is diagonalizable with all eigenvalues in the imaginary axis.
This is equivalent to the existence of a $s \times s$ invertible matrix $Q$ such that 
\begin{equation}
\label{eq:DD}
Q^{-1}\bar A Q = 
\left(
\begin{matrix}
0 & D \\
-D^T & 0 
\end{matrix}
\right)
\end{equation}
where $D$ is a real diagonal matrix (with non-negative diagonal entries) of size $m \times (s-m)$, where $m=[(s+1)/2]$ (and $s-m = [s/2]$). 

We will next show that (\ref{eq:DD}) may be exploited to solve efficiently linear systems of the form (\ref{eq:linsys}).  Consider the implicit equations (\ref{eq:Y2})--(\ref{eq:z}). Application of simplified Newton iteration  to such implicit equations leads to linear systems of the form 
\begin{equation}
\label{eq:linsys2}
 \begin{split}
\left( I_s \otimes I_\dim  - h \, \bar A \otimes J \right) \Delta Y - \half (e_s \otimes I_\dim) \, \Delta z &= r, \\
\left(  - h \, e_s^T B \otimes J \right) \Delta Y + \Delta z &= 0,
\end{split}
\end{equation}
where $e_s=(1,\dots,1)^T \in \R^{s}$. Clearly, if $(\Delta Y, \Delta z)$ is solution of (\ref{eq:linsys2}), then $\Delta Y$ is solution of (\ref{eq:linsys}).

By virtue of (\ref{eq:DD}), the linear system (\ref{eq:linsys2}) is equivalent, with the change of variables
$\Delta Y = (Q \otimes I_\dim) W$ to
\begin{equation}
\label{eq:linsys3}
\begin{split}
  \left(
  \begin{matrix}
    I_m \otimes I_\dim \ & \  -h \, D \otimes J  \\
    h \, D^T \otimes J \ & \ I_{s-m} \otimes I_\dim 
  \end{matrix}
\right) W -  \half\, (Q^{-1} e_s \otimes I_\dim) \, \Delta z &=  (Q^{-1} \otimes I_\dim )\, r, \\
- h \, (e_s^T  B Q \otimes J) \, W + \Delta z &= 0,
\end{split}
\end{equation}
The blockwise sparsity pattern of the system (\ref{eq:linsys3}) allows obtaining its LU decomposition by computing, in addition to several multiplications of matrices of size $\dim \times \dim$, the LU decompositions of $[s/2]+1$ real matrices of size $\dim \times \dim$: the matrices 
\[
I_\dim + h^2 \sigma_i^2 J^2, \quad i=1,\ldots,[s/2],
\]
 where $\sigma_1,\ldots,\sigma_{[s/2]} \geq 0$ are the diagonal entries in $D$, and an additional $\dim \times \dim$ matrix obtained from the former. We will give more details in  Subsection~\ref{ss:ssIRK} in the case of symmetric symplectic IRK schemes.

It is worth remarking that such a technique for solving linear systems of the form (\ref{eq:linsys}) is not restricted to symplectic IRK schemes. It is enough that the corresponding matrix $\bar A$  be diagonalizable with all its eigenvalues in the imaginary axis. This seems to be the case of several families of (non-symplectic) symmetric IRK methods of collocation type,  in particular,  for the nodes of Lobatto quadrature formulas, or if the nodes are either the zeros or the extrema of Chebyshev polynomials of the first kind.

\subsection{The case of symmetric symplectic IRK schemes}
\label{ss:ssIRK}

In the present section, in addition to the symplecticity conditions, that guarantee that the matrix $B^{\frac12}\bar A B^{-\frac12}$ is antisymmetric, we assume that the symmetry conditions (\ref{eq:symmetry_cond2}) hold. 

Consider the $s\times s$ orthogonal matrix $P=(P_1 \ P_2)$
such that, for $x =(x_1,\ldots,x_s)^T \in \R^s$, 
 $P_1^T x = (y_1 \cdots y_{m})^T$,   and $P_2^T x = (y_{m+1}, \cdots y_{s})^T$, 
where
\begin{align*}
y_i = \frac{\sqrt{2}}{2} (x_{s+1-i}+x_{i}), \quad &\mbox{for} \quad i=1,\ldots,[s/2], \\
y_{m} = x_{m}, \quad &\mbox{if} \quad s \mbox{ is odd}, \\
y_i = \frac{\sqrt{2}}{2} (x_{s+1-i}-x_{i}), \quad &\mbox{for} \quad i=m+1,\ldots,s,
\end{align*}
with $m=[(s+1)/2]$.

The symmetry condition (\ref{eq:symmetry_cond2}) implies that $P_i^T B^{\frac12} \bar A B^{-\frac12} P_i = 0$ for $i=1,2$, and since by symplecticity
$B^{\frac12} \bar A B^{-\frac12}$ is an antisymmetric matrix, we conclude that the matrix $\bar A$ is similar to
\begin{equation}
\label{eq:KK}
P^T B^{\frac12} \bar A B^{-\frac12} P = 
\left(
\begin{matrix}
0 & K \\
-K^T & 0 
\end{matrix}
\right)
\end{equation}
where $K = P_1^T B^{\frac12} \bar A B^{-\frac12} P_2$ (which is a real matrix of size $m \times (s-m) =[(s+1)/2]\times [s/2]$).
Let $K=U D V^T$ be the singular value decomposition of $K$, (where $U \in \R^{m \times m}$ and $V\in \R^{(s-m) \times (s-m)}$ are orthonormal matrices, and $D \in \R^{m \times (s-m)}$ is a diagonal matrix with the singular values $\sigma_1,\ldots,\sigma_{s-m}$ of $K$ as diagonal entries). We have that (\ref{eq:DD}) holds with  
\begin{equation*}
  Q = (Q_1 \ Q_2) = B^{-1/2}(P_1 \ P_2) \left(
\begin{matrix}
U & 0 \\
0 & V
\end{matrix}
\right) =  B^{-1/2} \left(
\begin{matrix}
 P_1 U &  P_2 V 
\end{matrix}
\right),
\end{equation*}
and $Q^{-1} = Q^T B$.
This implies that the linear system (\ref{eq:linsys2}), with the change of variables
\begin{equation}
\label{eq:DeltaYChVar}
 \Delta Y = (Q \otimes I_\dim) W = ( Q_1 \otimes I_\dim) W' + (Q_2 \otimes I_\dim) W''
\end{equation}
is equivalent to (\ref{eq:linsys3}). Due to the first symmetry conditions in (\ref{eq:symmetry_cond2}), $e_s^T B P_2 = 0$, and hence $e_s^T B Q_2 = e_s^T B P_2 V= 0$, so that (\ref{eq:linsys3}) reads
\begin{equation*}
  \begin{split}
       W' -h \, (D \otimes J)\,  W'' - \half\, (Q_1^T B e_s \otimes I_\dim) \, \Delta z &= (Q_1^T B \otimes I_\dim )\, r, \\
    h \, (D^T \otimes J)\, W'  + W'' \phantom{+\half\, (Q_2^T B e_s \otimes I_\dim) \, \Delta z  }
    &= (Q_2^T B \otimes I_\dim )\, r,\\
-h \, (e_s^T  B  Q_1 \otimes J) \, W' + \Delta z &= 0.
  \end{split}
\end{equation*}
 By solving for $W''$ from the second equation of the linear system above, 
\begin{equation}
\label{eq:W'}
 W'' = -h \, (D^T \otimes J)\, W'  + (Q_2^T B \otimes I_\dim )\, r.
\end{equation}
and substitution in the remaining two equations, one obtains
\begin{equation*}
\begin{split}
(I_{m} \otimes I_{\dim} + h^2\, D D^T \otimes J^2)  W' - \half\, ( Q_1^T B e_s \otimes I_\dim )\, \Delta z &= R, \\
-h \, (e_s^T  B Q_1 \otimes J) \, W' + \Delta z &= 0.
\end{split}
\end{equation*}
where $R=(Q_1^T B \otimes I_\dim) \, r + h \,  ( D Q_2^T B \otimes J)\,  r \in \R^{m\dim}$.

The linear system above can be rewritten in terms of 
\begin{equation*}
  R = \left(\begin{matrix}
R_1 \\ \vdots \\ R_m 
\end{matrix}\right), \quad   
W' = \left(\begin{matrix}
W_1 \\ \vdots \\ W_m 
\end{matrix}\right)
\end{equation*}
with $R_i, W_i \in \R^\dim$, $i=1,\ldots,m$, as follows:
\begin{align}
\label{eq:linsys4a}
(I_\dim +  h^2 \sigma_i^2 J^2) \, W_i - \frac{\alpha_i}{2}\,  \, \Delta z &= R_i, \quad i=1,\ldots,m,\\
\label{eq:linsys4b}
-h \, J \sum_{i=1}^{m} \alpha_i\,  W_i + \Delta z & = 0.
\end{align}
where 
\begin{equation*}
 \left(\begin{matrix}
\alpha_1 \\ \vdots \\ \alpha_m 
\end{matrix}\right) =  Q_1^T B e_s , 
\end{equation*}
and $\sigma_1 \geq \cdots \geq \sigma_{[s/2]}$ are the singular values of $K$, and  if $s$ is odd (in which case $m=[(s+1)/2]=[s/2]+1$),  then $\sigma_m=0$.

Thus, the unknown $\Delta z \in \R^{\dim}$ can be obtained by solving the linear system
\begin{equation}
\label{eq:linsysz}
  M \,  \Delta z = h\,  J \sum_{i=1}^{m} \alpha_i  (I_\dim +  h^2 \sigma_i^2 J^2)^{-1} R_i,
\end{equation}
where 
\begin{equation}
\label{eq:M}
M = I_\dim + J\, \frac{h}{2} \, \sum_{i=1}^m \alpha_i^2 (I_\dim +  h^2 \sigma_i^2 J^2)^{-1}   \in \R^{\dim \times \dim}.
\end{equation}

The unknowns in $W' \in \R^{m\dim}$ are then solved from (\ref{eq:linsys4a}), while $W'' \in \R^{(s-m)\dim}$ may be obtained from (\ref{eq:W'}).

The required solution $\Delta Y$ of the original linear system (\ref{eq:linsys}), may finally be obtained from (\ref{eq:DeltaYChVar}).

\subsection{Alternative reformulation of symplectic IRK schemes}
\label{ss:alternative}

If the coefficients  $b_i,a_{i j}$ determining a symplectic IRK are replaced by floating point numbers $\tilde b_i,\tilde a_{i j}$ that approximate them, 
then the resulting IRK scheme typically fails to satisfy the symplecticity conditions (\ref{eq:sympl_cond_1}). This results~\cite{Hairer2008} in a  method that exhibits a linear drift in the value of quadratic invariants of the system and in the Hamiltonian function when applied to autonomous Hamiltonian systems.

Motivated by that,  the  map $\Phi:\R^{d+2} \to \R^{d}$  of the one-step integrator  (\ref{eq:one_step})  corresponding to the IRK scheme, defined by (\ref{eq:PhiIRK})--(\ref{eq:Y}), is rewritten in~\cite{Antonana2017} in the following equivalent form:
\begin{equation*}
\label{eq:PhiIRK2}
\Phi(y,t,h) :=y + \sum_{i=1}^s L_{i},
\end{equation*}   
where $L_{i} \in \R^d$, $i=1,\ldots,s$ are implicitly defined as functions of $(t,y,h) \in \R^{d+2}$ by 
\begin{equation}
\label{eq:L}
 L_{i} = h \, b_i \, f(t+c_i h, y+ \sum_{j=1}^s \mu_{ij}\,L_{j}), \quad  i=1 ,\ldots, s, 
 \end{equation}
where
\begin{equation*} 
\mu_{ij}=a_{ij}/b_j,  \quad 1 \leq i,j \leq s.
\end{equation*}
The symplecticity condition (\ref{eq:sympl_cond_1}) is equivalent to
\begin{equation} \label{eq:sympl_cond_2}
\mu_{ij}+\mu_{ji}-1=0, \quad 1 \leq i,j \leq s.
\end{equation}

The main advantage of the proposed formulation over the standard one is that the absence of multiplications in the symplecticity condition (\ref{eq:sympl_cond_2}) makes possible to find machine number approximations $\mu_{i j}$ of $a_{i j}/b_j$ satisfying exactly the symplecticity condition (\ref{eq:sympl_cond_2}).

With that alternative formulation, the Newton iteration reads as follows:  Initialize $L_i^{[0]}=0$ for $i=1,\ldots,s$, and compute for $k=1,2,\ldots$
\begin{equation}
\label{eq:Newton_iteration2}
\begin{split}
1) 
   & \quad Y_i^{[k]} := y+\sum_{j=1}^{s}\, \mu_{ij}\, L_{j}^{[k-1]}, \quad  i=1 ,\ldots, s. \\
   & \quad g_i^{[k]} := -L_{i}^{[k-1]}  + h \, b_i\, f(t+c_i h,  Y_i^{[k]}), \quad  i=1 ,\ldots, s, \\   
2) & \quad \mathrm{Solve \ } \Delta L_{i}^{[k]} \mbox{\  from \ } \\
   & \quad \Delta L_{i}^{[k]}  - h b_i J_i^{[k]} \sum_{j=1}^{s}\, \mu_{ij} \, \Delta L_{j}^{[k]}=g_i^{[k]}, \quad  i=1 ,\ldots, s,  \\
& \mbox{where} \quad  J_i^{[k]}=\frac{\partial f}{\partial y}(t + c_i h,Y_{i}^{[k]}) \quad	\mbox{for} \quad  i=1,\ldots, s,  \\
3) 
   & \quad   L^{[k]} := L^{[k-1]}  + \Delta L^{[k]},
   \end{split}
\end{equation}
In the simplified version of the Newton iteration where the Jacobian matrices $J_i^{[k]}$ are replaced by a common approximation $J$ (say, $J=\frac{\partial f}{\partial y}(t + h/2, y)$), the linear system in (\ref{eq:Newton_iteration2}) is replaced by
\begin{equation}
\label{eq:Newton_iteration_sls}
\Delta L^{[k]} = \left( I_s \otimes I_\dim  - h \, B A B^{-1} \otimes J \right)^{-1} 
\left(
\begin{matrix}
g_1^{[k]}\\
\vdots\\
g_s^{[k]}
\end{matrix}
\right),
\end{equation}
In that case, we need to repeatedly solve systems of the form 
\begin{equation}
\label{eq:linsysZG}
\left( I_s \otimes I_\dim  - h \, B A B^{-1} \otimes J \right) \, \Delta L = g, 
\end{equation}
for prescribed $g \in \R^{s\dim}$. Repeated solution of linear systems of this form is also required if the linear system in (\ref{eq:Newton_iteration2}) is iteratively solved as described in Subsection~\ref{ss:iterativeLS} below.

Of course, (\ref{eq:linsysZG}) can be solved by adapting the technique described in Subsections~\ref{ss:sIRK} and \ref{ss:ssIRK}  for the solutions of systems of the form (\ref{eq:linsys}).
We next describe, for the symmetric case (i.e., when the symmetry condition (\ref{eq:symmetry_cond}) holds),  the corresponding procedure (with the notation adopted in Subsection~\ref{ss:ssIRK}) to compute the solution $\Delta L$ of  (\ref{eq:linsysZG}):
\begin{enumerate}
\item LU decompositions:
\begin{itemize}
\item Compute the LU decompositions of the $\R^{\dim\times \dim}$ matrices 
\begin{equation}
\label{eq:IsigmaJ2}
I_\dim+ h^2\, \sigma_i^2\, J^2, \quad i=1,\ldots,[s/2],
\end{equation}

\item  Compute the matrix  $M \in \R^{\dim\times \dim}$ given in (\ref{eq:M}) (recall that $\sigma_m=0$ when $s$ is odd), and obtain its LU decomposition.
\end{itemize}

\item Solution of system (\ref{eq:linsysZG}):
\begin{itemize}
 \item Compute $R \in \R^{m\dim}$ from 
\begin{equation*}
R=(Q_1^T  \otimes I_\dim) \, g + h \,  ( D Q_2^T \otimes J)\, g,
\end{equation*}
\item Compute
\begin{equation*}
d= h\, J \sum_{i=1}^{m} \alpha_i  (I_\dim +  h^2 \sigma_i^2 J^2)^{-1} R_i,
\end{equation*}
\item Compute $\Delta z \in \R^{\dim}$ as the solution of the linear system $M\, \Delta z = d$,
\item  Next, compute $W_1,\ldots,W_m \in \R^{\dim}$ from 
\begin{equation*}
(I_\dim +  h^2 \sigma_i^2 J^2) \, W_i - \frac{\alpha_i}{2}\, J \, \Delta z = R_i, \quad i=1,\ldots,m.
\end{equation*}
\item Follow by computing $W_{m+1},\ldots,W_s \in \R^{\dim}$  from 
\begin{equation*}
 \left(
\begin{matrix}
W_{m+1}\\
\vdots\\
W_s
\end{matrix}
\right)
  = - \left(
\begin{matrix}
h \, \sigma_1 \, J \, W_1\\
\vdots\\
h \, \sigma_{s-m}\, J\, W_{s-m} 
\end{matrix}
\right)
+ (Q_2^T\otimes I_\dim )\, g.
\end{equation*}
\item And finally, $\Delta L \in \R^{s\dim}$ is obtained from 
\begin{equation*}
 \Delta L = (B Q \otimes I_\dim)
  \left(
\begin{matrix}
W_{1}\\
\vdots\\
W_s
\end{matrix}
\right).
\end{equation*}
\end{itemize}
\end{enumerate}

\section{Implementation of symplectic IRK schemes with Newton-like iterations}

In this section, we present an algorithm that implements symplectic IRK schemes by making use of the techniques in previous section. Special care is taken to try to reduce the effect of round-off errors by adapting some techniques used in~\cite{Antonana2017} for the implementation of symplectic IRK schemes with fixed point iterations.
Our algorithm is intended to be applied with the 64-bit IEEE double precision floating point arithmetic.

\subsection{Auxiliary techniques}
\label{ss:auxiliary}

In this subsection we summarize some techniques associated to the use of finite precision arithmetic that we applied in the fixed point iteration implementation of symplectic IRK schemes proposed in~\cite{Antonana2017}, and  will be used in the algorithm proposed in Subsections~\ref{ss:implementation}.

Let $\F \subset \R$ 
 be the set of machine numbers of the 64-bit IEEE double precision  floating point arithmetic. We consider the map $\fl:\R \longrightarrow \F$   that sends each real number $x$ to a nearest machine number $\fl(x) \in \F$ .

\subsubsection{Kahan's compensated summation}
\label{sss:Kahan}

The application of any one-step integrator of the form (\ref{eq:one_step}) requires computing sums of the form
\begin{equation}
\label{eq:sumy_n}
y_{n+1} = y_{n} + x_n, \quad n=0,1,2,\ldots, 
\end{equation}
For an actual implementation that only uses a floating point arithmetic with machine numbers in $\F$, special care must be taken with the additions (\ref{eq:sumy_n}). The naive recursive algorithm 
$\hat y_{n+1} :=\fl(\hat y_{n} + \fl(x_n))$,  ($n=0,1,2,3\ldots$),
typically suffers, for large $n$,  a significant loss of precision due to round-off errors.  It is well known that such a round-off error accumulation can be greatly reduced with the use of Kahan's compensated summation algorithm~\cite{Kahan1965} (see also~\cite{Higham2002}\cite{Muller2009}). 

 Given $y_0 \in \R^{\dim}$ and a sequence $\{x_0,x_1,\ldots,x_n,\ldots\} \subset \F^{\dim}$ of machine numbers, Kahan's algorithm is aimed to compute the sums  $y_n = y_0 + \sum_{\ell=0}^{n-1} x_{\ell}$, ($n\geq 1$,) 
  using a prescribed floating point arithmetic, more precisely than with the naive recursive algorithm.
  The actual algorithm reads as follows: 
  
  \begin{algorithm}[H]
   \BlankLine
    $\tilde y_0= \fl(y_0); \ e_0=\fl(y_0-\tilde y_0)$\;
    \BlankLine
    \For{$l\leftarrow 0$ \KwTo $n$}
    {
     \BlankLine
      $X_l = \fl(x_l + e_{l})$\;
      $\tilde y_{l+1} = \fl(\tilde y_{l} + X_{l})$\;
      $\hat X_{l} = \fl(\tilde y_{l+1} - \tilde y_{l})$\;
      $e_{l+1} =  \fl(X_{l} - \hat X_{l})$\;
     \BlankLine
    }
   \caption{Kahan’s compensated summation}
   \label{alg:Kahan'sCS}
  \end{algorithm}
 The sums $\tilde y_l + e_l \in \R^{\dim}$  are more precise approximations of the exact sums $y_l$ than $\tilde y_l \in \F$. 
 Algorithm \ref{alg:Kahan'sCS} can be interpreted as a family of maps parametrized by $n$ and $\dim$,
 \begin{equation*}
S_{n,\dim}:\F^{(n+3)\dim} \to \F^{2\dim},
\end{equation*}
  that given the arguments $\tilde y_0,e_0,x_0,x_1,\ldots,x_n \in \F^\dim$,  returns $\tilde y_{n+1}, e_{n+1} \in \F^\dim$ such that $\tilde y_{n+1} + e_{n+1} \approx (\tilde y_0 + e_0) + x_0+ x_1 + \cdots + x_n$ with some small error.

\subsubsection{Stopping criterion for iterative processes}
\label{sss:stopping}

Given a  smooth map $F:\R^D \to \R^D$ and $Z^{[0]} =(Z_1^{[0]},\ldots,Z_D^{[0]}) \in \R^D$ assume that the iteration
\begin{equation}
\label{eq:Ziter}
Z^{[k] }=F(Z^{[k-1]}), \quad \mbox{for} \quad k=1,2,\ldots
\end{equation}
produces a sequence  $\{ Z^{[0]},Z^{[1]}, Z^{[2]}, \ldots \} \subset \R^D$ that converges to a fixed point $Z^{[\infty]}$ of $F$.

Assume now that instead of the original map $F$, we have a computational substitute 
\begin{equation}
\label{eq:tildef}
\widetilde F: \F^D\to \F^D.
\end{equation}
Ideally, for each $Z \in \F^D$, $\widetilde F(Z):= \fl(F(Z))$. In practice, the intermediate computations to evaluate $\widetilde F$ are typically made using the floating point arithmetic corresponding to $\F$, which will result in some additional error  caused by the accumulated effect of several round-off errors. 

The resulting sequence $\widetilde{Z}^{[k] }=\widetilde{F}(\widetilde{Z}^{[k-1]})$, $k=1,2,\ldots$ (started with  $\widetilde{Z}^{[0]} = \fl(Z^{[0]})$) will either converge to a fixed point of $\widetilde{F}$ in a finite number $K$ of iterations or will fail to converge. In the former case, the fixed point $\widetilde{Z}^{[K]} \in \F^D$ of $\widetilde{F}$ may be expected to be a good approximation of the fixed point $Z^{[\infty]} \in \R^D$ of $F$. In the later case, one would expect that there exists an index $K$ such that the approximations $\widetilde{Z}^{[k]}\approx Z^{[\infty]}$ improves for increasing $k$ up to $k=K$, but the quality of the approximations $\widetilde{Z}^{[k]}$ does not improve for $k>K$. It then make sense to apply an algorithm of the form

 \begin{algorithm}[H]
  \BlankLine
   $ k=0$\;
   $ \widetilde{Z}^{[0]}=\fl(Z^{[0]})$\;
   \BlankLine
   \While{ $( \ \mathrm{ContFcn}(\widetilde Z^{[0]}, \cdots, \widetilde Z^{[k]}) \ )$ }
   {
    \BlankLine
     $k= k+1$\;
     $\widetilde{Z}^{[k]}=\widetilde{F}(\widetilde{Z}^{[k-1]})$\;
    \BlankLine
   }
  \caption{Stopping criterion}
  \label{alg:ZiterK}
 \end{algorithm}
where $\mathrm{ContFcn}(\widetilde{Z}^{[0]}, \cdots, \widetilde{Z}^{[k]}
)$  gives either true if it is judged appropriate to continue iterating, and false otherwise.
In ~\cite{Antonana2017}, we propose defining this function so that  $\mathrm{ContFcn}(\widetilde{Z}^{[0]}, \cdots, \widetilde{Z}^{[k]}
)$ returns
\begin{align*}
         \left\{ 
            \begin{array}{lcl}
            \text{false} &\rightarrow & \ \mathrm{if} \quad (\widetilde Z^{[k]} = \widetilde Z^{[k-1]}) \ \text{or} \\
                         & & \quad \quad  k>1 \text \ {and} \ k=K-1 \ {and} \ k=K \text \ {and} \  \forall  j \in \{1,\ldots,D\},\\
                         & & \quad \quad \min \left(\{| \widetilde{Z}_j^{[1]}-\widetilde{Z}_j^{[0]}|,\cdots ,| \widetilde{Z}_j^{[k-1]}-\widetilde{Z}_j^{[k-2]}|\} \ /\{0\} \right) \leq | \widetilde{Z}_j^{[k]}-\widetilde{Z}_j^{[k-1]}|\\
            \\
            \text{true} &\rightarrow & \ \mathrm{otherwise}.  
            \end{array} 
         \right. 
 \end{align*}

The output $\widetilde{Z}^{[K]}$ of the algorithm will  be a fixed point of $\widetilde{F}$  when it stops because  $\widetilde Z^{[K]}= \widetilde Z^{[K-1]}$, and in any case it is  not expected that $\widetilde{Z}^{[k]}$ for $k>K$ be a better approximation of the fixed point $Z^{[\infty]} \in \R^D$ of $F$ than $\widetilde{Z}^{[K]}$.

\subsection{An inexact Newton iteration}
\label{ss:iterativeLS} 

In our implementation of symplectic IRK schemes to be described in Subsection~\ref{ss:implementation}, we consider a modified version of the Newton iteration (\ref{eq:Newton_iteration2}). In each iteration, instead of exactly solving for $\Delta L_{i}^{[k]}$ 
from a linear system of the form
\begin{equation}
\label{eq:linsysZ}
 \Delta L_{i}^{[k]}  - h b_i \ J_i \sum_{j=1}^{s}\, \mu_{ij} \, \Delta L_{j}^{[k]} = g_i^{[k]}, \quad  i=1 ,\ldots, s, 
\end{equation}
where
\begin{equation}
\label{eq:gk}
g_i^{[k]} = -L_{i}^{[k-1]}  + h \, b_i\, f\Big(t+c_i h,  y+\sum_{j=1}^{s}\, \mu_{ij}\, L_{j}^{[k-1]}\Big), \quad  i=1 ,\ldots, s, \\
\end{equation}
and
\begin{equation*}
\Delta L^{[k]}= \left(
\begin{matrix}
\Delta L_1^{[k]}\\
\vdots\\
\Delta L_s^{[k]}
\end{matrix}
\right) \in \R^{s\dim}, \ \
g^{[k]}= \left(
\begin{matrix}
g_1^{[k]} \\
\vdots\\
g_s^{[k]}
\end{matrix}
\right) \in \R^{s\dim},
\end{equation*}
we iteratively compute a sequence $ \Delta L_{i}^{[k,0]},  \Delta L_{i}^{[k,1]},  \Delta L_{i}^{[k,2]}, \ldots$ of approximation of its solution $\Delta L^{[k]} \in \R^{s\dim}$ as shown below (Algorithm \ref{alg:InexactNw}).

\begin{algorithm}[H]
  $ \Delta L^{[k,0]} = (I_s \otimes I_d - h \ BAB^{-1} \otimes J)^{-1} \ g^{[k]}$\;
  \BlankLine
   \While{ $( \ \mathrm{ContFcn}(\fl_{32}(\Delta L^{[k,0]}), \cdots,\fl_{32}(\Delta L^{[k,\ell]})) \ ) $ }
  {
   \BlankLine
   $l=l+1$\;
   $G_i^{[k,\ell]} = g_i^{[k]} - \Delta L_{i}^{[k,\ell-1]}  + h b_i \ J_i \sum_{j=1}^{s}\, \mu_{ij} \, \Delta L_{j}^{[k,\ell-1]},  \quad  i=1 ,\ldots, s$\;
   \BlankLine
   $\Delta L^{[k,l]}=\Delta L^{[k,l-1]}+ (I_s \otimes I_d - h \ BAB^{-1} \otimes J)^{-1} \ G^{[k,l]}$\;
  }
 \caption{Inner iteration}
 \label{alg:InexactNw}
\end{algorithm}
Hereafter, $\fl_{32}(x)$ denotes the 32-bit IEEE single precision machine number that is closest to $x \in \R$, and we let $\fl_{32}$ act componentwise on vectors.

\paragraph*{}In the algorithm we propose in Subsection~\ref{ss:implementation},  the Jacobian matrices $J_i$ in (\ref{eq:linsysZ}) will be evaluated in approximations of the stage values $Y_i$ that  are accurate at the single precision level. This implies that it does not make sense to apply the iteration (Algorithm \ref{alg:InexactNw}) until an accurate double precision approximation $\Delta L^{[k,\ell]}$ of the solution  $\Delta L^{[k]}$ of (\ref{eq:linsysZ}) is obtained. Motivated by that, 
we will stop the iteration when 
$\mathrm{ContFcn}(\fl_{32}(\Delta L^{[k,0]}), \cdots,\fl_{32}(\Delta L^{[k,\ell]}))$ returns false (i.e., typically, when $\fl_{32}(\Delta L^{[k,\ell]}) =\fl_{32}(\Delta L^{[k,\ell-1]})$.

\subsection{Algorithm for one step of the IRK scheme}
\label{ss:implementation}

In our implementation, the numerical solution $y_n \approx y(t+h n) \in \R^{\dim}$, $n=1,2,\ldots$, is obtained as the sum  $\tilde y_n + e_n $ of two vectors in $\F^{\dim}$.  In particular,  the initial value $y_0 \in \R^\dim$ is (approximately) represented as $\tilde y_0 + e_0$, where  $\tilde y_0 = \fl(y_0)$ and $e_0=\fl(y_n-\tilde y_n)$. 
Instead of (\ref{eq:one_step}), we actually have 
\begin{equation*}
(\tilde y_{n+1},e_{n+1}) = \tilde \Phi(\tilde y_n, e_n, t_n, t_{n+1}-t_n),
\end{equation*}
where $\tilde \Phi: \F^{2\dim+2} \to \F^{2\dim}$. 

Our proposed  implementation of one step 
\begin{equation*}
(\tilde y^*, e^*) = \tilde \Phi(\tilde y,e,t,h)
\end{equation*}
of the IRK scheme is performed in five substeps:
 \begin{enumerate}
\item 
Starting from $L^{[0]}=0 \in \R^{s\dim}$, we apply several simplified Newton iterations (i.e., the simplified version of Newton iterations (\ref{eq:Newton_iteration2}) where the linear system is replaced by (\ref{eq:Newton_iteration_sls})) to compute 
\[
L^{[1]}=L^{[0]} +\Delta L^{[1]}, \quad  L^{[2]}=L^{[1]} +\Delta L^{[2]} , \ldots
\]
 until $\fl_{32}(L^{[k]}) =\fl_{32}(L^{[k-1]})$ (or rather, by using the notation introduced in paragraph~\ref{sss:stopping},
until  $\mathrm{ContFcn}(\fl_{32}(L^{[0]}), \cdots,\fl_{32}(L^{[k]}))$ returns false).

\item
Use $L^{[k]}$ to compute the Jacobian matrices 
\begin{equation*}
J_i= \frac{\partial f}{\partial y}\left(t +c_i\,  h,  \ \tilde y+  \sum_{j=1}^{s}\, \mu_{ij}\, L_{j}^{[k]}\right), \quad i=1,\ldots,s,
\end{equation*}
\item 
Then consider the increment $\Delta L^{[k]} \in \F^{s\dim}$ obtained in first substep as an approximation $\Delta L^{[k,0]}$ of the exact solution $\Delta L^{[k]}$ of the linear system in  (\ref{eq:Newton_iteration2}), and 
apply the inner iterations (Algorithm \ref{alg:InexactNw}) to obtain as output an approximation $\Delta L^{[k,\ell]}$ (accurate at least at single precision level).
\item 
Follow by updating $L^{[k]} = L^{[k-1]} + \Delta L^{[k,\ell]}$, and $k=k+1$, and
applying a final inexact Newton iteration with the Jacobian matrices $J_i$ computed in the second substep. More preciselly,
compute an approximation $\Delta L^{[k,\ell]}$ (again accurate at least at single precision level) of 
the solution $\Delta L^{[k]}$ of  (\ref{eq:linsysZ})--(\ref{eq:gk}) by applying Algorithm~\ref{alg:InexactNw}.
\item
Finally, the increment $\tilde \Phi(\tilde y,e,t,h)$,  defined as the sum $(\tilde y + e) + \sum_{i=1}^{s}(L_i^{[k-1]} + \Delta L^{[k,\ell]})$ is accurately obtained as the (unevaluated) sum of the double precision vectors $\tilde y^*,e^*\in \F^\dim$ with the help of Kahan's  compensated summation algorithm (summarized in paragraph~\ref{sss:Kahan}) as follows: First, perform the sum $\delta :=  e + \sum_{i=1}^{s} \Delta L^{[k,\ell]})$ of the vectors with relatively smaller size in the double precision floating point arithmetic, and then compute   $(\tilde y^*,e^*) = S_{s,d}(\tilde y, \delta, L_1^{[k-1]}, \ldots, L_s^{[k-1]})$.
\end{enumerate}

Some remarks about our actual implementation are in order:
\begin{itemize}
\item All the linear system with coefficient matrix $(I_s \otimes I_\dim  - h \, B A B^{-1} \otimes J )$ are solved by means of the algorithm at the end of Section~3.
\item The coefficients $\mu_{i j}$ are machine numbers in $\F$ (i.e., in the target precision floating point system) satisfying exactly the symplecticity condition (\ref{eq:sympl_cond_2}) and the symmetry conditions $\mu_{j, i} = \mu_{s+1-i,s+1-j}$.
\item  The remainders  (\ref{eq:gk}) ($i=1,\ldots,s$, $k\geq 1$) should in principle be computed with $y\in \R^\dim$ replaced by  $\tilde y + e$ ($\tilde y, e \in \F^\dim$). However, the effect of ignoring the extra digits of $y$ that may be contained in $e$ is expected to be so small that it should be enough to take it into account only in the final inexact Newton iteration (substep 4 above). That is, it should be enough considering (\ref{eq:gk}) with $y\in \R^\dim$ replaced by $\tilde y\in \F^d$ in all the Newton-like iterations with the exception of the final one. And in the final inexact Newton iteration,  rather than computing  (\ref{eq:gk}) with $y\in \R^\dim$ replaced by $\tilde y + e$, we make use of the Jacobian matrices $J_i$ to obtain the following approximation
\begin{equation*}
 h \, b_i\, f\Big(t+c_i h,  \tilde y + e +\sum_{j=1}^{s}\, \mu_{ij}\, L_{j}^{[k-1]}\Big)  -L_{i}^{[k-1]} \approx 
\left(h \, b_i\, f_i^{[k]} -L_{i}^{[k-1]} \right) + h\, b_i \, J_i \, e,
\end{equation*}
where $f_i^{[k]}=f\Big(t+c_i h,  \tilde y +\sum_{j=1}^{s}\, \mu_{ij}\, L_{j}^{[k-1]}\Big)$.

\item  If the FMA (fused-multiply-add) instruction is available, it should be used to compute $h \, b_i\, f_i^{[k]} -L_{i}^{[k-1]}$ (with precomputed coefficients $h  b_i \in \F$ satisfying the symmetry conditions $h  b_{s+1-i}  = h  b_i  $).
\end{itemize}

Our final implementation
is summarized in Algorithm~\ref{alg:step}.  

 \begin{algorithm}[h]
  \BlankLine
    $L^{[0]} = 0$; \   $J= \frac{\partial f}{\partial y}(t +h/2,\tilde y)$\;
    $M = I_\dim + J\, \frac{h}{2} \, \sum_{i=1}^m \alpha_i^2 (I_\dim +  h^2 \sigma_i^2 J^2)^{-1}$\;
    $\text{Compute the LU decomposition of M}$\;
    \BlankLine
    $\text{/************************** 1st substep **************************/}$\;
    \BlankLine
    $k=0$\;
   \While{ $\mathrm{ContFcn}(\fl_{32}(L^{[0]}), \ldots, \fl_{32}(L^{[k]}))$}
   {
    \BlankLine
    $k=k+1$\;        
    $Y_{i}^{[k]} = \tilde y+\sum_{j=1}^{s}\, \mu_{ij}\, L_{j}^{[k-1]}, \ i=1,\dots,s$\; 
    $f_{i}^{[k]} =  f\left(t+c_i h, Y_{i}^{[k]} \right), \ i=1,\dots,s$\;
    $g_i^{[k]} =h \, b_i\,f_i^{[k]} -L_i^{[k-1]}, \ i=1,\dots,s$\;
    $\Delta L^{[k]} =\left(I_s \otimes I_\dim  - h \, B A B^{-1} \otimes J \right)^{-1} g^{[k]}$\;
    $L^{[k]} = L^{[k-1]} + \Delta L^{[k]}$\;    
    \BlankLine
   }
   \BlankLine
    $\text{/************************** 2nd substep **************************/}$\;
   \BlankLine 
    $J_i= \frac{\partial f}{\partial y}\left(t + c_i h, \ \tilde y+\sum_{j=1}^{s}\, \mu_{ij}\, L_{j}^{[k]}\right), \ i=1,\dots,s$\;
   \BlankLine 
    $\text{/************************** 3rd substep **************************/}$\;
   \BlankLine
    $ \ell=0$\;
    $\Delta L^{[k,0]} =  \Delta L^{[k]}$\;
    \While{$ \mathrm{ContFcn}(\fl_{32}(\Delta L^{[k,0]}), \ldots, \fl_{32}(\Delta L^{[k,\ell]}))$}
    {
    \BlankLine
    $\ell=\ell+1$\;
    $ G_i^{[k,\ell]} = g_i^{[k]} - \Delta L_{i}^{[k,\ell-1]}  + h b_i J_i \sum_{j=1}^{s}\, \mu_{ij} \, \Delta L_{j}^{[k,\ell-1]}, \ i=1,\dots,s$\;
    $\Delta L^{[k,\ell]} =\Delta L^{[k,\ell-1]} + \left(
         I_s \otimes I_\dim  - h \, B  A B^{-1} \otimes J \right)^{-1} 
          G^{[k,\ell]}$\;       
     \BlankLine
   }
    $L^{[k]} = L^{[k-1]} + \Delta L^{[k,\ell]}$\;    
    \BlankLine 
    $\text{/************************** 4th substep **************************/}$\;
    $k=k+1$\;
    $Y_{i}^{[k]} = \tilde y+ \sum_{j=1}^{s}\, \mu_{ij}\, L_{j}^{[k-1]}, \ i=1,\dots,s$\;
    $f_{i}^{[k]} =  f\left(t+c_i h, Y_{i}^{[k]} \right), \ i=1,\dots,s$\;
    $g_i^{[k]} = \left(h \, b_i\,f_i^{[k]} -L_i^{[k-1]} \right) + h \, b_i \, J_i \, e, \ i=1,\dots,s$\;
    $ \ell=0$\;
    $\Delta L^{[k,0]} =\left(
    I_s \otimes I_\dim  - h \, B A B^{-1} \otimes J \right)^{-1} g^{[k]}$\;
    \While{$ \mathrm{ContFcn}(\fl_{32}(\Delta L^{[k,0]}), \ldots, \fl_{32}(\Delta L^{[k,\ell]}))$}
    {
    \BlankLine
    $\ell=\ell+1$\;
    $ G_i^{[k,\ell]} = g_i^{[k]} - \Delta L_{i}^{[k,\ell-1]}  + h b_i J_i \sum_{j=1}^{s}\, \mu_{ij} \, \Delta L_{j}^{[k,\ell-1]}, \ i=1,\dots,s$\;
    $\Delta L^{[k,\ell]} =\Delta L^{[k,\ell-1]} + \left(
         I_s \otimes I_\dim  - h \, B A B^{-1} \otimes J \right)^{-1} 
          G^{[k,\ell]}$\;       
    \BlankLine
    }
   \BlankLine 
    $\text{/************************** 5th substep **************************/}$\;
    \BlankLine
     $\delta = e + \sum_{i=1}^{s} \Delta L_i^{[k,\ell]}$\;
     $(\tilde y^*,e^*) = S_{s,d}(\tilde y, \delta, L_1^{[k-1]}, \ldots, L_s^{[k-1]})$\;
    \BlankLine
  \caption{Implementation of one step of the IRK scheme}
  \label{alg:step}
 \end{algorithm}

\clearpage

\section{Numerical experiments}
\label{s:ne}

We next report some numerical experiments to asses our implementation of the $6$-stage Gauss collocation method of order $12$ based on Newton-like iterations (Algorithm~\ref{alg:step})  with 64-bit IEEE double precision  floating point arithmetic. 
\subsection{The double pendulum stiff problem}
We consider the planar double pendulum stiff problem: a double bob pendulum with masses $m_1$ and $m_2$ attached by rigid massless rods of lengths $l_1$ and $l_2$ and spring of elastic constant $k$ between both rods (the rods are aligned at equilibrium). For $k=0$, the problem is non-stiff, and the system's stiffness arises through increasing the value of $k$.

The configuration of the pendulum is described by two angles $q=(\phi,\theta)$ : while $\phi$ is the angle (with respect to the vertical direction) of the first bob with, the second bob's angle is defined by $\psi=\phi+\theta$. We denote the corresponding momenta as $p=(p_{\phi},p_{\theta})$.

Its Hamiltonian function $H(q,p)$ is
\begin{multline} \label{eq:2}
-\frac{ {l_1}^2 \ (m_1+m_2) \ {p_{\theta}}^2 +{l_2}^2 \ m_2 \ (p_{\theta} -p_{\phi})^2 + 2 \ l_1 \ l_2 \ m_2 \ p_{\theta} \ (p_{\theta} -p_{\phi}) \  \cos(\theta )} {{l_1}^2  \ {l_2}^2 \ m_2 \  (-2 \ m_1 - m_2 + m_2 \ \cos(2 \theta ))} \\
-g  \ \cos (\phi) \  (l_1 \ (m_1+m_2)+l_2 \ m_2 \ \cos(\theta))+g \ l_2 \ m_2 \ \sin(\theta) \sin(\phi)+\frac{k}{2} \ \theta^2 ,
\end{multline}
We consider the following fixed parameters values
\begin{equation*} \label{eq:17}
g=9.8 \ \frac{m}{sec^2}\ ,\ \ l_1=1.0 \ m \ , \ l_2=1.0 \ m\ , \ m_1=1.0 \ kg\ , \ m_2=1.0 \ kg,
\end{equation*} 
and we choose different values for the elastic constant $k$ to study different levels of stiffness in the double pendulum. 
The initial values are chosen as follows: For $k=0$, we choose the initial values considered in~\cite{Antonana2017} that gives rise to a non-chaotic trajectory,  $q(0)=(1.1, -1.1)$ and $p(0)=(2.7746,2.7746)$. The initial values for $k\neq 0$ are chosen as
\begin{equation*}
q(0)=\left(1.1, \frac{-1.1}{\sqrt{1+100k}}\right), \quad
p(0)=(2.7746,2.7746)
\end{equation*}
so that the total energy of the system is bounded as $k \to \infty$.

All the integrations are performed with step-size $h=2^{-7}$, that is small enough for round-off errors to dominate over truncation errors in the non-stiff case $k=0$. The truncation errors dominate over the round-off errors for large enough stiffness constant $k>0$.
We have integrated over $T_{end}=2^{12}$ seconds and sample the numerical results every $m=2^{10}$ steps. 

\subsection{Round-off error propagation}

First, we check the good performance of round-off error propagation of our new implementation based on Newton-like iterations. In \cite{Antonana2017} we proposed an implementation based on fixed-point iterations for non-stiff problems that takes special care of reducing the propagation of round-off errors. We will compare the round-off error of both implementations of the $6$-stage Gauss collocation method.

We have studied in detail the errors in energy of the double pendulum problem for three values of $k$: $k=0$, where the round-off errors dominate over truncation errors, $k=2^{10}$, where both kinds of errors are similar in size,  and  $k=2^{12}$, where truncation errors dominate over round-off errors. In order to make a more robust comparison of the numerical errors due to round-off errors, we adopt  (as in~\cite{Hairer2008}) an statistical approach.  We have considered for each of the three initial value problems, $P=1000$ perturbed initial values by randomly perturbing each component of the initial values with a relative error of size $\mathcal{O}(10^{-6})$.

\begin{figure}[h!]
\centering
\begin{tabular}{c c}
\subfloat[$k=0$: mean energy error ]
{\includegraphics[width=.4\textwidth]{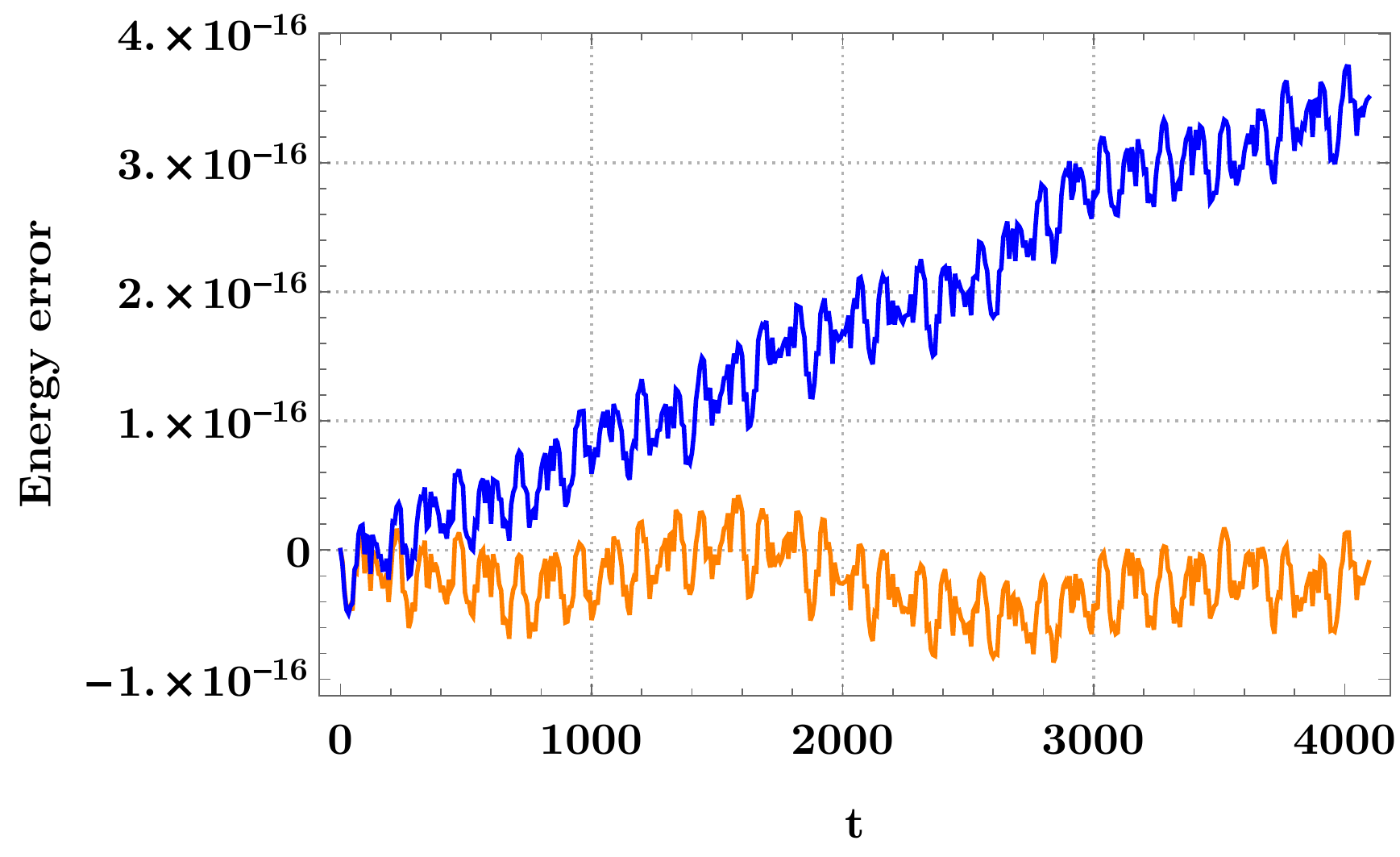}}
&
\subfloat[$k=0$: standard  deviation energy error]
{\includegraphics[width=.4\textwidth]{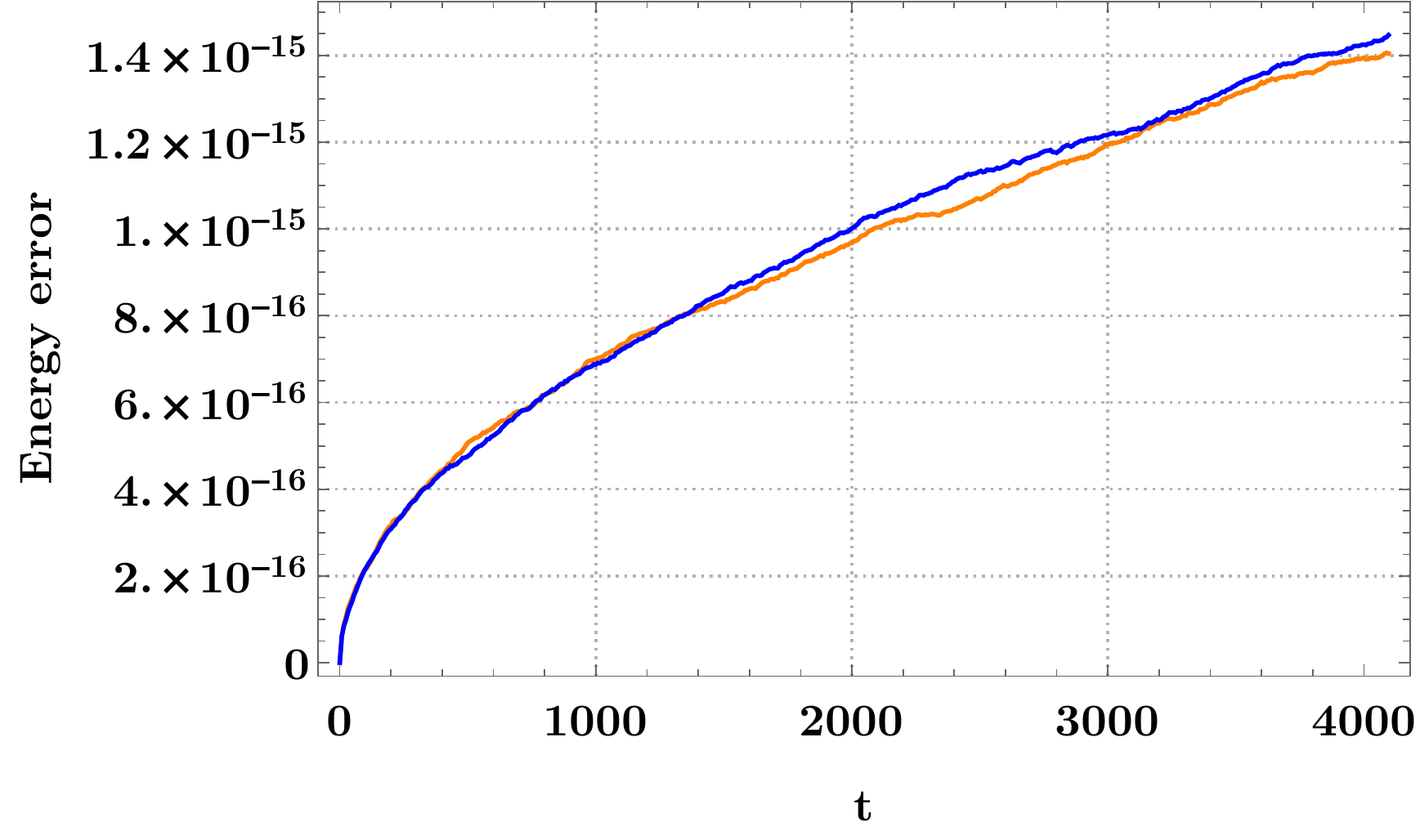}}
\\
\subfloat[$k=2^{10}$: mean energy error ]
{\includegraphics[width=.4\textwidth]{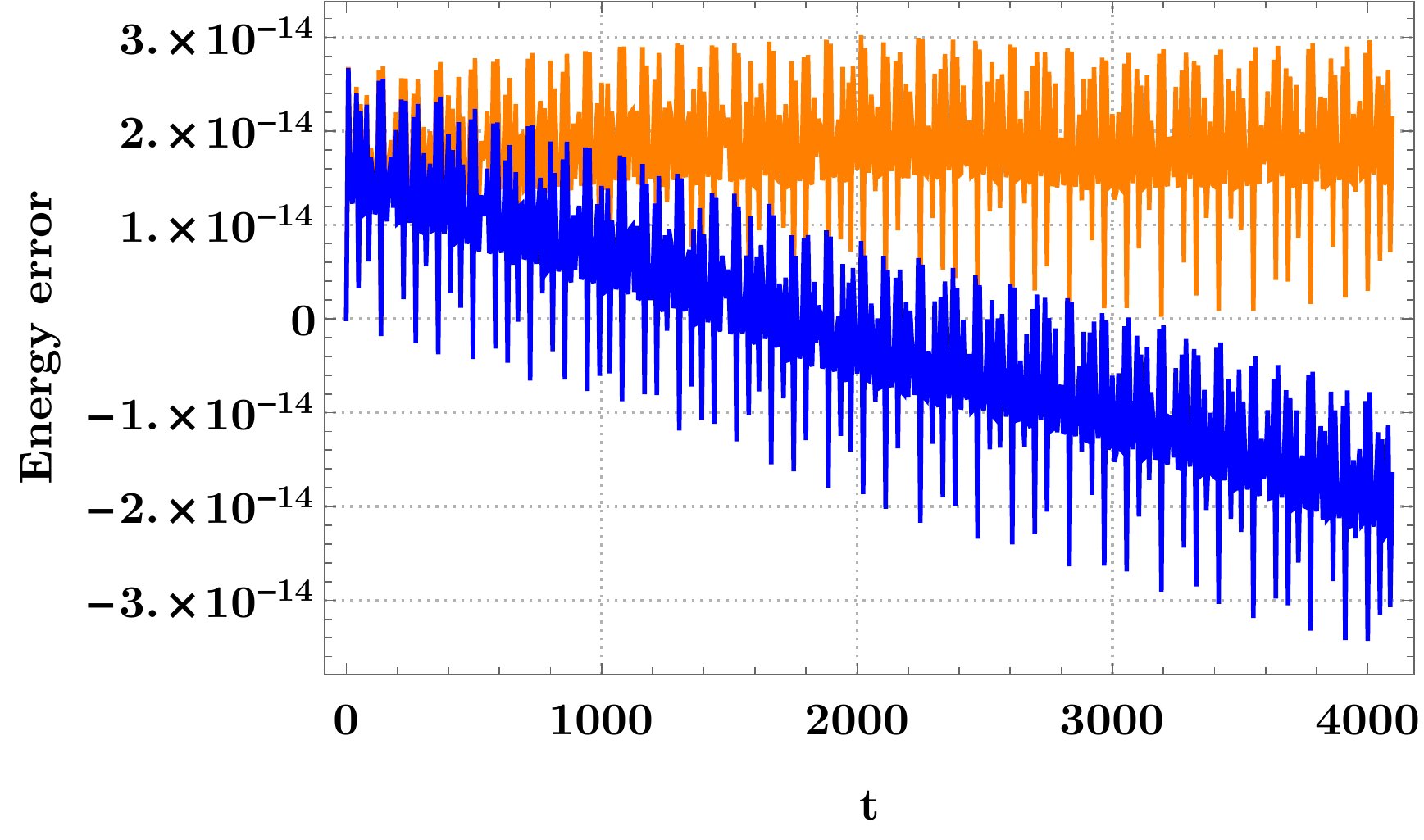}}
&
\subfloat[$k=2^{10}$: standard deviation energy error]
{\includegraphics[width=.4\textwidth]{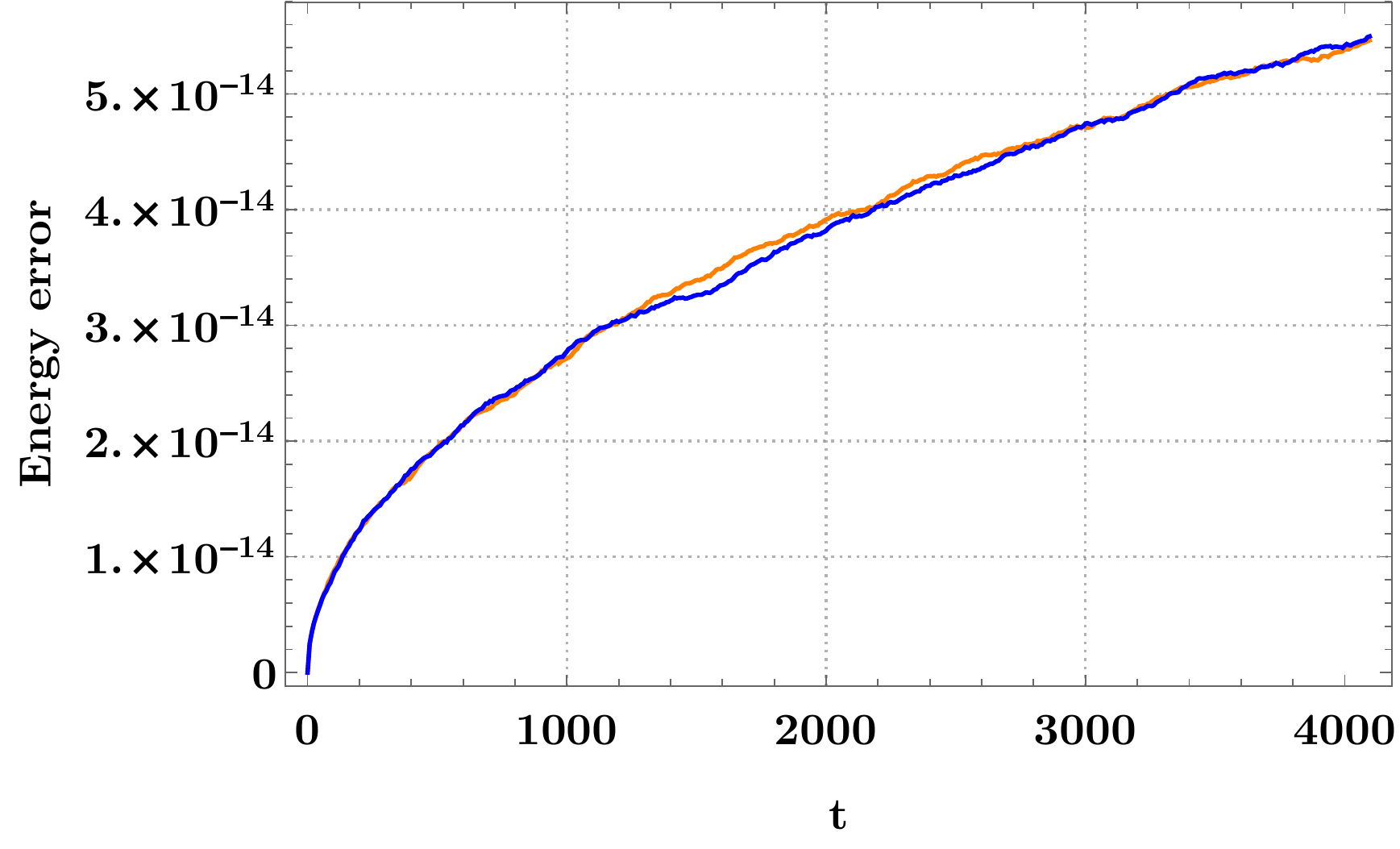}}
\\
\subfloat[$k=2^{12}$: mean energy error]
{\includegraphics[width=.4\textwidth]{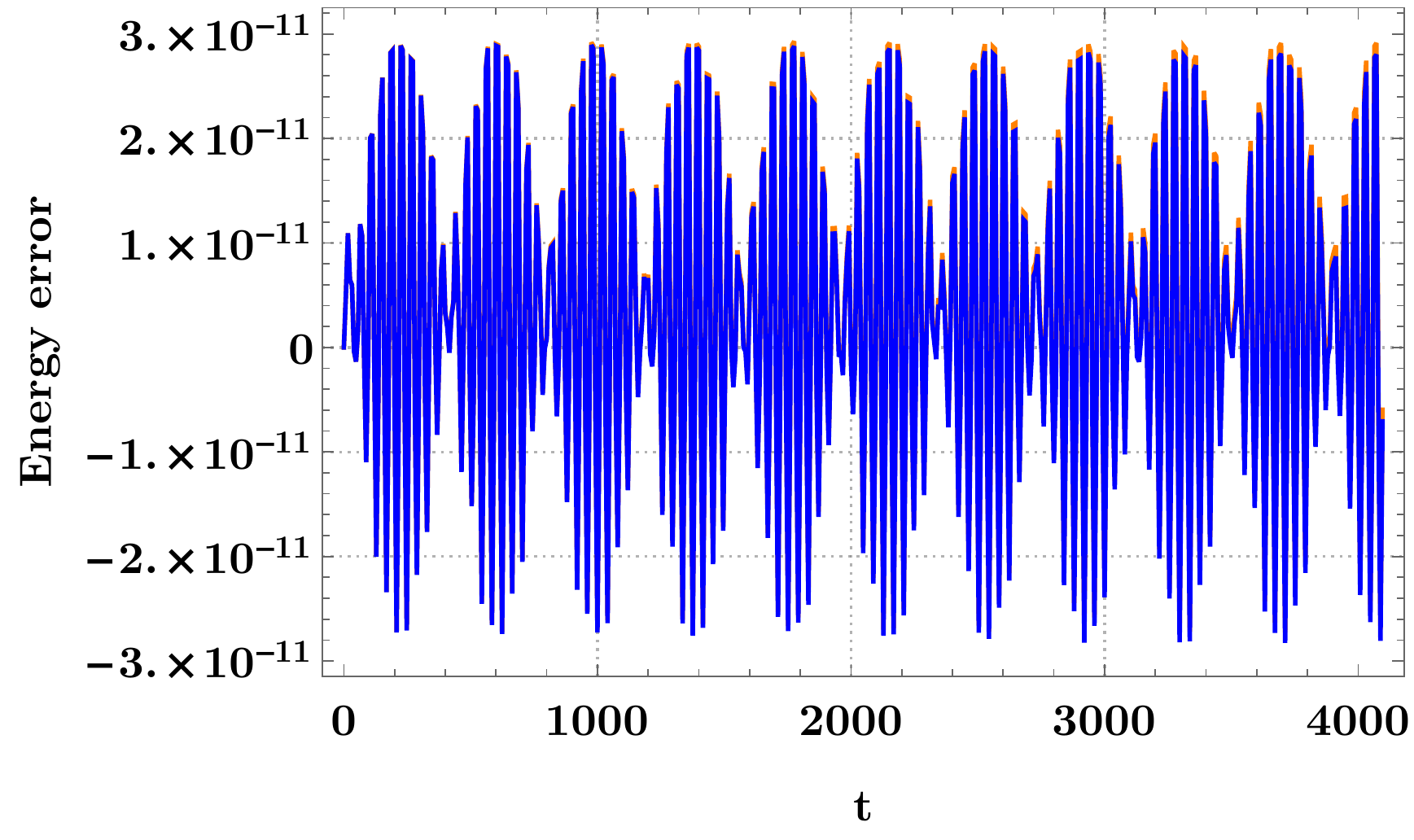}}
&
\subfloat[$k=2^{12}$: standard deviation energy error]
{\includegraphics[width=.4\textwidth]{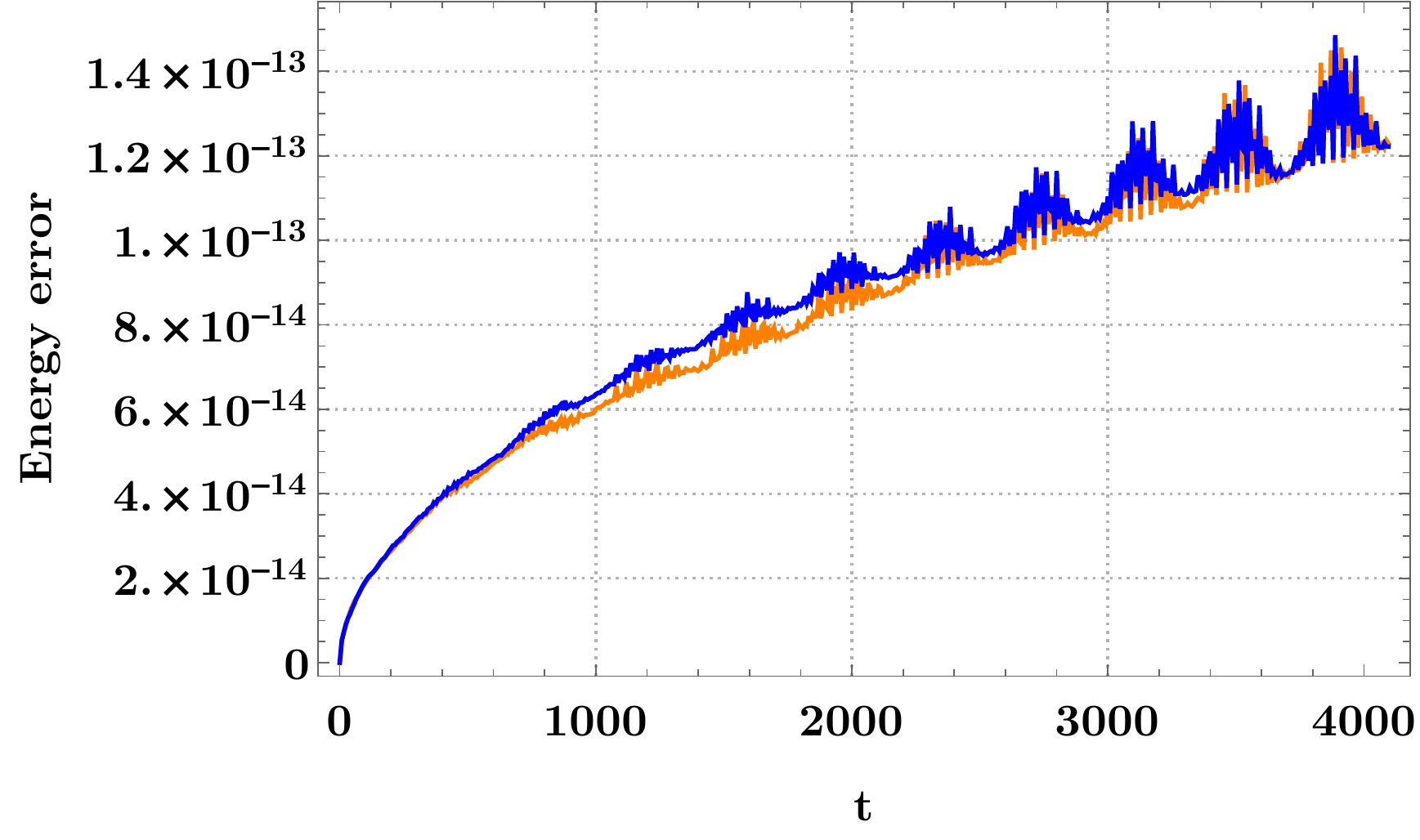}}
\end{tabular}
\caption{\small Evolution of mean  (left) and standard deviation (right) of relative errors in energy for fixed-point implementation (blue), and  Newton implementation (orange). Non-stiff case $k=0$ ~ (a,b), one stiff case $k=2^{10}$ ~(c,d) and second stiff case $k=2^{12}$ ~(e,f)}
\label{fig:plot1}
\end{figure}

The numerical tests in Figure~\ref{fig:plot1} seem to confirm the good performance of round-off error propagation of our new implementation. From one hand, one can observe that, as in~\cite{Antonana2017}, the fixed-point implementation exhibits a small linear drift of the mean energy error for $k=0$ and $k=2^{10}$, while in the Newton implementation this energy drift does not appear at all. On the other hand, the standard deviation of the energy errors are of similar size and grow proportionally to $t^{\frac12}$ in both implementations. 

\subsection{Fixed-point versus Newton iteration}

We summarize in Table~\ref{tab:fp1} the main results of numerical integrations for both implementations:
the fixed-point iteration and Newton-like iteration for four different values of $k$.

We have compared their efficiency by sequential execution of each iteration method, and reported the cpu-time of each numerical integration. 
In addition, we have reported  the number of iterations per step (It. per step) in both implementations and the number of linear systems solved in the Newton implementation. To check the precision of the numerical solution, we have reported the maximum relative energy error,
\begin{equation*}
 \max \left| \frac{E(t_i)-E(t_0) }{E(t_0)} \right|, \quad t_i=t_0+ i h, \ \ i=0,\dots,steps.
\end{equation*}

\paragraph*{}We can see that for low  values of $k$, the fixed-point implementation is more efficient than Newton implementation. But as we increase the stiffness of the double pendulum, the number of iteration needed at each step in the fixed-point implementation grows up notably, while in the Newton implementation the number of iterations even becomes slightly lower for higher values of $k$. Hence, the Newton implementation eventually becomes more efficient as the stiffness increases. For $k$ values higher than $k=2^{18}$, the fixed-point iteration fails to converge, while the Newton implementations succeeds while keeping approximately the same number of iterations per step (cpu-time: 17s; iterations per step: 4.95; Linear solves per step: 10.94).

\begin{table}[h!]
\caption[Summary of numerical integrations] 
{\small{Summary of numerical integrations with fixed-point iteration and Newton iteration based implementation for the following spring's elastic values: $k=0, \ k=2^6 , \ k=2^{12}$ and $\ k=2^{16}$. $E_0$ indicates the initial energy of the system.  We show the cpu-time, the number of iteration per step (It. per step), the number of linear system solving operations (L. solves per step) and maximum energy error for each numerical computation}}
\label{tab:fp1}       
\centering
{%
\begin{tabular}{ lcccc } 
 \hline
\\
 k               & $0$  & $2^6$ & $2^{12}$ & $2^{16}$ \\
 $E_0$           & $-14.39$  & $-5.75$ & $-5.64$ & $-5.64$ \\ 
\\
 \hline
\\
 Fixed-points it.&           &         &         &         \\
 \cline{1-1}     &           &         &         &         \\
 Cpu-time (sec.)    & $10$      & $12$    & $19$    & $51$    \\ 
 It. per step    & $8.58$    & $11.1$  & $22.$  & $64.2$  \\
 Max energy error  & $2.96\times 10^{-15}$ & $1.81\times 10^{-14}$ & $2.94\times 10^{-11}$ & $6.33\times 10^{-5}$ \\
 \\
 Newton it.            &           &         &         &         \\
 \cline{1-1}           &           &         &         &         \\
 Cpu-time (sec.)   & $18$      & $20$    & $19$    & $18$     \\
 It. per step          & $5.09$    & $5.53$  & $5.58$  & $5.01$   \\
 L. solves per step    & $11.37$   & $12.92$  & $12.72$  & $11.04$ \\
 Max energy error      & $1.6\times 10^{-15}$ & $1.74\times 10^{-14}$ & $2.94\times 10^{-11}$ & $6.33\times 10^{-5}$ \\   
 \\  
   \hline
 \end{tabular}}
\end{table}

\section{Conclusions}

Our main contribution is a technique to solve efficiently the simplified linear systems of symplectic IRK schemes. This technique can be adapted for some symmetric non-symplectic schemes as well. Such technique could also be exploited for the numerical solution of boundary value problems with collocation methods with Gaussian quadrature nodes.

In addition, an efficient algorithm for implementing symplectic IRK methods with reduced round-off error propagation is provided.
A C-code with our implementation for $s$-stage Gauss collocation method of order $2s$ in the 64-bit IEEE double precision floating point arithmetic can be downloaded from IRK-Newton Github software repository or go to the next url: \url{https://github.com/mikelehu/IRK-Newton}.

\paragraph{Acknowledgements}
M. Anto\~nana, J. Makazaga, and A. Murua have received funding from the Project of the Spanish Ministry of Economy and Competitiveness with reference MTM2016-76329-R (AEI/FEDER, EU), from the
project MTM2013-46553-C3-2-P from Spanish Ministry of Economy and Trade, 
and as part of the Consolidated Research Group IT649-13 by the Basque Government.

\end{document}